\documentclass[reqno]{amsart}
\usepackage{amssymb}

\newtheorem{theorem}{Theorem}[section]
\newtheorem{oldtheorem}{Theorem}

\newtheorem{lemma}[theorem]{Lemma}

\theoremstyle{definition}
\newtheorem{definition}[theorem]{Definition}
\newtheorem{problem}[theorem]{Problem}

\theoremstyle{remark}
\newtheorem{remark}[theorem]{Remark}
\newtheorem{corollary}[theorem]{Corollary}

\theoremstyle{plain}
\newtheorem{conjecture}[theorem]{Conjecture}

\numberwithin{equation}{section}

\newcommand{\bvp}{boundary value problem }
\newcommand{\bvps}{boundary value problems }
\newcommand{\Grf}{Green's function }
\newcommand{\efet}{entire function of exponential type}
\newcommand{\UB}{\textup{(UB) }}
\newcommand{\UBP}{\textup{(UBP) }}
\newcommand{\fss}{fundamental system of solutions}
\newcommand{\mcm}{modified characteristic matrix}
\DeclareMathOperator{\entier}{entier}

\providecommand{\abs}[1]{\lvert#1\rvert}
\providecommand{\norm}[1]{\lVert#1\rVert}

\newcommand{\mbb}[1]{\mathbb{#1}}
\newcommand{\Cb}{\mbb{C}}
\newcommand{\Rb}{\mbb{R}}
\newcommand{\mbf}[1]{\mathbf{#1}}
\newcommand{\A}{\mbf{A}}
\newcommand{\B}{\mbf{B}}
\newcommand{\F}{\mbf{F}}
\newcommand{\Q}{\mbf{Q}}
\newcommand{\D}{\mbf{D}}
\newcommand{\gD}{\mbf{\Delta}}
\newcommand{\gT}{\mbf{\Theta}}
\newcommand{\gP}{\mbf{\Psi}}

\newcommand{\bl}{\pmb{[}}
\newcommand{\br}{\pmb{]}}

\begin{document}

\title[SPECTRALITY of OPERATORS]{Spectrality of ordinary differential operators}

\author[A.Minkin]{Arkadi Minkin}
\address{A.Minkin, str. Haarava 8/5, 36863 Nesher, Israel}
\email{arkadi\_minkin@yahoo.com}

\dedicatory{To my parents Sarah and Moisei}

\begin{abstract}
We prove the long standing conjecture 
in the theory of two-point boundary value problems
that completeness and Dunford's spectrality imply Birkhoff regularity. 
In addition we establish the even order part of S.G.Krein's conjecture that
dissipative differential operators are Birkhoff-regular and give sharp estimate 
of the norms of spectral projectors in the odd case.

Considerations are based on a new direct method, 
exploiting \textit{almost orthogonality} of Birkhoff's solutions of the equation $l(y)=\lambda y$. 
 This property was discovered earlier  by the author.
\end{abstract}

\maketitle

\tableofcontents

\section{Introduction}
\subsection{Dunford's spectrality.}
Question of unconditional convergence of spectral decompositions
plays a central role in the spectral theory. It is also known under other frameworks: 
similarity of a linear operator to a normal one,
N.Dunford's spectrality, 
free interpolation problem 
\cite[appendix A.3]{NV98},
\cite{DunSchwIII,NikOper2,NikShift}.
In this paper we deal with complete operators with compact resolvent. For them the questions above
translate respectively into unconditional basicity \UB of eigenfunctions (ef), of root subspaces
and of eigen and associated functions (eaf).
The latter is the most general setting. So we shall refer to it as the \UB problem, though sometimes will use instead the term spectrality as is common in the theory of \bvps (briefly bvps).

Set $D=-id/dx$ and consider 
a \bvp (bvp) in $L^2(0,1)$,
defined by a differential expression 
\begin{equation}
l(y)\equiv D^ny+\sum_{k=0}^{n-2}p_k(x)D^ky=\lambda y,\ \ 0\le x\le 1,\ \
p_k\in L(0,1)  \label{eq1.1a}
\end{equation}
and $n$ linearly independent boundary conditions 
\begin{equation}
U_{\jmath }(y)\equiv \sum_{k=0}^{n-1} 
  \left( 
    a_{jk}D^jy(0) + b_{jk}D^jy(1)
  \right) =0, 
  \;\; j=0,\ldots ,n-1.  \label{eq1.1b}
\end{equation}
\noindent
Spectral theory of the operator $L$, defined by this bvp, is thoroughly investigated during the last hundred years.
The bibliography is enormous and we shall refer the reader to the fundamental
monographs of M.A.Naimark \cite{NaiLdoEngI,NaiLdoEngII} and N.Dunford,J.T.Schwartz \cite{DunSchwIII}.

Remind that inverse of $L$ is a finite-dimensional perturbation of a Volterra operator
\begin{equation}
\label{eq-volt-pert}
   L^{-1}f=Vf+\sum_{j=0}^{n-1}(f,h_j)g_j
\end{equation}
where $Vf$ gives solution to the Cauchy problem for $l(y)$ 
and zero boundary conditions
\begin{equation}
   \label{eq-zbc}
   D^jy(0)=0, \ j=0,\ldots,n-1. 
\end{equation}
Spectral investigation of this class of linear operators was initiated by A.P.Hromov \cite{hro74}.
During the last 20 years question of similarity to normal for
operators (\ref{eq-volt-pert}) with dissipative $V$
was deeply explored in works of G.M.Gubreev \cite{gub00,gub03}  and we will dwell on it further. 

The aim of this paper is three-fold. 
First, to give a final solution to the \UB problem for two-point bvps \eqref{eq1.1a}-\eqref{eq1.1b}. 
Second, explain why their spectrality don't fit into all existing schemes.
For the reader's convenience we will show this on simplest examples.
As a byproduct we give an account of relevant abstract results as well as those obtained in the classical spectral theory of bvps.

And the last but not least, we 
believe that solution of the \UB problem for (\ref{eq1.1a})-(\ref{eq1.1b}) will help
in investigation of much more difficult class (\ref{eq-volt-pert}) with non-dissipative Volterra operator $V$.
\subsection{Paper outline.}
In the section \ref{sec:exp} we describe results on basisness of exponentials and ef
 of first and higher order bvps.
 It is continued by an account of the Stone-regularity in the section \ref{sec:stone} and abstract approaches in the section \ref{sec:abstract}.
We consider all methods from the viewpoint of bvps, trying to reveal obstacles that
prevent their usage for solving the \UB problem.
Add that modern \emph{projection method} grew from the theory of exponentials.
Obviously the latter is still important as a source of ideas, technique and inspiration for further investigations.
 
Background for spectral theory of bvps
is exposed in the section \ref{sec:green}.
Note that it includes a new notion of a
\mcm{ } and a new definition of regularity determinants.
 
Section \ref{sec:main} describes main results of the paper together with some 
open conjectures. In addition,
we derive partial solution of the S.G.Krein's conjecture
about spectrality of dissipative differential operators
directly from one of our main results, theorem \ref{thm-2}.

Proofs are placed in sections \ref{sec:mcm}-\ref{sec:proof-thm1}. 
In the section \ref{sec:mcm} we establish properties of the \mcm{.}
Here we deduce theorem \ref{thm-2}, which ties together minimal resolvent's 
growth with nonvanishing of regularity determinants.
Then in the section \ref{sec:growth} we establish density properties
of eigenvalues, lying in a sector off the real axis.
These results form a foundation of the proof of the main theorem \ref{thm-1}
in the section \ref{sec:proof-thm1}.
\subsection{Notations.}
Throughout the paper components of matrices and vectors are enumerated
beginning from zero.
Matrices are written in boldface together with their brackets
to distinguish such bracket from Birkhoff's symbol, e.g.
\[
\gD=\bl\Delta_{jk}\br _{j,k=0}^{n-1}.
\]
Different constants are denoted $C, C_1, c$ and so on. They may vary even
during a single computation.
Other notations and abbreviations:
\begin{itemize}
\item  $[a]:=a+O(1/\varrho )$ stands for the Birkhoff's symbol;
\item fss - \fss;
\item $\Cb_\pm$ - upper/lower half-plane, $\Rb$ - real axis;
\item $H_\pm^2$  - Hardy space in $\Cb_\pm$;
\item ev - eigenvalue(s);
\item cv - characteristic value(s);
\item efet - \efet{;} 
\item $\left| \Delta \mathop{\longleftarrow}\limits_{k} d\right|$
stands for determinant $\Delta$ with the $k$-th
column replaced by a vector $d$;
\item  $A\asymp B$ means a double-sided estimate $\ \ C_1\cdot |A|\leq
|B|\leq C_2\cdot |A|\ \ $ with some absolute constants $C_{1,2}$, which
don't depend on the variables $A$ and $B$.
\end{itemize}
\section{Exponentials and bvps.}
\label{sec:exp}
\subsection{First order bvp.}
\label{subsec:1st-ord}
Note that for $l(y)\equiv{D}y$ and general functional $U(y)$ in the boundary condition we arrive at the classical question of unconditional basicity of exponentials. 
The first, now classical results in this direction 
were initiated  by Paley and Wiener \cite{PW}.
They discovered that harmonic frequencies $k$ in the orthogonal
system $\{\exp(ikx)\}_{-\infty}^\infty$ may be replaced by \textit{close} real ones $\lambda_k$,
preserving $\{\exp(i\lambda_k x)\}_{-\infty}^\infty \in (UB)$ in
 $L^2(0,2\pi)$.

At the beginning of 1960ies B.Ya.Levin and V.D.Golovin established basis properties 
of exponentials whose generating function is of sine type \cite{le61,go63}. 
In particular, this method implies unconditional basisness with parentheses \UBP of the first-order bvp
\begin{align} 
		Ly   &=  Dy	
		\label{eq0.1} \\
		U(y) &=  \int_{0}^{a} y(x) d\sigma(x)	    \label{eq0.2} 		
\end{align}
with finite measure $d\sigma$ on $[0,1]$, provided that
\begin{equation}
	\sigma\{0\} \ne 0, \qquad \sigma\{a\} \ne 0.	\label{eq0.3}
\end{equation}
The ef are exponentials 
$e_{kj}=x^j\exp(i\lambda_kx)$,
$\lambda_k$ are ev, and we meet here an example of a bvp with rather general functional in the boundary condition.
Afterward the notion of a sine-type function led  to 
one of descriptions of Riesz bases 
from exponentials \cite{aj88}. Remind that Riesz basis is an unconditional almost normalized basis,
i.e. $C_1 \le \| e_{kj}\|_{L^2(0,a)}\le C_2$. Therefore necessarily frequencies lie in a strip
$	|\Im\lambda_k|\le C $. This implicit requirement is removed when
one passes to unconditional bases from exponentials.
\subsubsection{Discovery of projection method.}
However, already in 1973 B.S.Pavlov devised a simple geometric approach to this problem \cite{pa73}, using functional model of dissipative operators. For differential operator in a finite interval his method  
requires that Fourier transform of a functional in the boundary condition
be a sine-type function. This requirement seems difficult to verify.

Therefore B.S.Pavlov performed a deep investigation of the asymptotic behavior of this Fourier transform 
and of the ev in case of the boundary condition (\ref{eq0.2}),
assuming piecewise-absolute continuity of the measure $d\sigma$. As a result he obtained an \textit{effective} criterion of Riesz basisness of the root vectors in their span 
 \cite[Theorem 5]{pa73}. It means that under appropriate conditions
 the carlesonity of the spectrum was \underline{proved} 
 and not merely put into
 theorem's assumptions.
\subsubsection{Semi-bounded spectrum.}
Afterward Pavlov's ideas grew to the \emph{projection method} \cite{pa79}, which solved the \UB problem 
for exponentials in a finite interval $[0,a]$,
provided that the spectrum $\Lambda=\{ \lambda_k\}_1^\infty$
is semi-bounded
\begin{equation}
	\delta = \inf \Im\Lambda > -\infty, 
	\quad \Im\Lambda:=\{ \Im\lambda, \lambda\in\Lambda \}.
	\label{eq-semi-bounda}
\end{equation}
Recall the criterion
\cite[end of p.658]{pa79}, modulo N.K.Nikolskii's remark 
cited therein at the end of p.658, see also
\cite{nik80}\footnote{In this paper N.K.Nikolskii extended projection method to systems from values of reproducing kernels.}.
Let 
\[
b_\lambda(z):= \frac{ \abs{\lambda^2+1} }{\lambda^2+1}
\frac{z-\lambda}{z-\overline{\lambda}},\; \lambda\ne i;
\quad
b_i(z) = \frac{z-i}{z+i}
\] 
be the Blaschke factor. Choose some $\eta >\delta$.
Then the criterion consists of three conditions: 
\begin{align}
	& \Lambda+i\eta \in (C),\; \text{i.e. }
	\inf_k
	\prod_{\substack{
	                  j=1 \\
	                  j\ne k}
	      }            	     
	     ^\infty 
	   \abs{ 
	   b_{\lambda_k+\eta}(\lambda_j+\eta)	  	   
	   } 
	   > 0.
	\label{eq_carl}  \\
	& \abs{\varphi(x+i\eta)}^2 \in (A_2), \label{eq_a2} \\
	& \varphi(z) \ \text{is efet with indicator diagram } [0,a]. \label{eq_a3}
\end{align}
Here $\varphi(z)$ denotes the \emph{generating function} of $\Lambda$, i.e. entire function with zeros only in $\Lambda$ counting multiplicities. 
The sum in $\Lambda+\eta$ is element-wise.
$(C)$ and $(A_2)$ stand for Carleson and Muckenhoupt conditions,
see definitions in \cite{garbaf}.
Obviously (\ref{eq-semi-bounda}) may be replaced by
\begin{equation}
	\inf \Im\Lambda > 0 \label{eq-semi-bound}
\end{equation}
via multiplication by $\exp(\eta x)$ for any fixed $\eta > \inf\Im\Lambda$.
\subsubsection{Arbitrary spectrum.}
A general criterion of \UB of exponentials without spectrum restriction has been obtained in \cite{min91}. The formulation is essentially the same as in 
(\ref{eq_carl})-(\ref{eq_a3}) but for a new set of frequencies, obtained by reflecting
$\lambda_k \in \Cb_{-}$ to the upper half-plane, while $\lambda_k \in \Cb_{+}\cup\Rb$ stay the same.
Later the criterion was transferred to the \UBP case in \cite{min94}. These works used
arithmetics of coinvariant subspaces of inverse shift. 
Note that \cite{min91} provides also another form of the criterion via distances 
of an appropriate unimodular symbol similar to N.K.Nikolskii's theorem \cite{nik80}, and thus covers incomplete systems of exponentials with arbitrary spectrum,
constituting \UB in the span. 

This situation appears in applications \cite{vlaiva}
or when studying
exponential bases in scales of interpolation spaces, for instance, in Sobolev spaces
\cite{aviva,ivakal}.
\subsubsection{Interpolating sequences.}
Apparently the \UB problem for exponentials may be restated as interpolation problem.
Namely, by Paley-Wiener theorem Fourier transform maps $L^2(0,a)$ to the classical Paley-Wiener space
$PW^2$ of efet with indicator diagram in $[0,a]$,
square summable on $\Rb$.
Thus we come to the question:
when does interpolation problem
$
    f(\lambda_k) = a_k
$
have a unique solution $f\in PW^2$ for every data $\{a_k \}$ satisfying
\[
  \sum_k \abs{a_k}^2 / \norm{ \exp(i\lambda_k x) }_{L^2(0,a)}^2 < \infty.
\]
Such set $\Lambda$ is called a \textit{complete interpolating sequence} (CIS).
K.Seip and Yu.I.Lyubarskii found another proof \cite{seilu97} of the main theorem from 
\cite[theorem 0.1, equivalence $A)\Leftrightarrow C)$]{min91}. 
They reduced the problem to the boundedness of a discrete Hilbert transform, which got the following remarkable solution:
\[
    \abs{\varphi(x)/dist(x,\Lambda)}^2\in (A_2).
\]
Here we don't need to reflect frequencies into $\Cb_{+}$.
However the case of incomplete sequences remained uncovered.
It seems that it is related to the essence of their method,
see the dichotomy conjecture in  \cite[p.717]{seiIcm98}.
\subsubsection{Criterion for functional-measure.}
Note that these theorems are stated in terms of spectrum distribution and of the behaviour of the generating function. It is quite natural from the function theory viewpoint, where \underline{the given data is $\Lambda$}.
At the same time for bvps the given data are \underline{boundary conditions},
and it is needed to establish a result directly in their terms.
For instance, the following theorem is valid.
\begin{oldtheorem}
\label{thm-A}
Let $d\sigma$ be a discrete measure on $[0,a]$.
Then for $(\ref{eq0.1})-(\ref{eq0.2}) \in (UBP)$
it is necessary and sufficient that (\ref{eq0.3})
be fulfilled.
\end{oldtheorem}
Sufficiency belongs to B.Ya.Levin-V.D.Golovin,
whereas necessity \cite{min97} was obtained, using
the general criterion from \cite{min94}.
Hence in \cite{min97} it was shown how to check carlesonity of the spectrum for 
such general bvp.
Actually (\ref{eq0.3}) is nothing 
else but its Birkhoff-regularity.
Note also that M.Rubnich extended theorem
\ref{thm-A} to a measure with discrete and  singular continuous components
\cite{ru96}.
\subsection{Higher order bvps.}
Now let us turn to ordinary differential operators.
With abuse of notations assume that boundary conditions (\ref{eq1.1b}) 
are normalized
\cite{sal68}:
\begin{equation}
U_{\jmath }(y)\equiv b_{\jmath }^0D^jy(0)+b_{\jmath }^1D^jy(1)+\ldots
=0,\;\;j=0,\ldots ,n-1.  \label{eq1.2}
\end{equation}
\noindent
The ellipsis takes place of lower order terms at $0$ and at $1$;
$b_{\jmath }^0,b_{\jmath }^1$ are column vectors of length $%
r_{\jmath }$, where
\begin{equation}
\label{eq-rankB}
  0\le r_{\jmath }\le 2,\;
  \sum_{k=0}^{n-1}r_{\jmath }=n,\ \
  rank\left(b_{\jmath }^0b_{\jmath }^1\right)=r_{\jmath }.
\end{equation}
Evidently $r_{\jmath }=0$ implies
absence of order $j$ conditions. In the case $r_{\jmath }=2$ we put
\[
\left( b_{\jmath }^0b_{\jmath }^1\right) =\left(
\begin{array}{cc}
1 & 0 \\
0 & 1
\end{array}
\right) .
\]
Below in the subsection \ref{std_reg} we define Birkhoff-regular boundary conditions. 
They possess a lot of remarkable spectral properties: estimate of
the Green's function, asymptotics of ev and ef, 
equconvergence with trigonometric Fourier series on any compact $[a,b]\subset(0,1)$.
Moreover recently there were established necessary and sufficient conditions for equiconvergence on the whole interval $[0,1]$, see \cite[$n=2$]{hro75} and 
\cite[Chapter 2, $n\ge 2$]{min99a}.

However, only in \textup{1960}ties G.M.Kesel\'{}man \cite{kes64} and V.P.Mihailov \cite{mikh62}
proved that \emph{strong regularity}
(briefly $L\in (SR)$), see definition \ref{s-reg},
yields \UB of eaf.
 
In definition \ref{s-reg} we call 
\emph{Birkhoff} but not \emph{strongly regular} boundary conditions
\emph{weakly regular} for evident reasons and write $L\in (WR)$.
For them A.A.Shkalikov established
unconditional basicity with parentheses (two summands in each) \cite[1979]{shk79a}.
Obviously in this case spectrality vanishes if eaf are not paired, see examples of
G.M.Kesel\'{}man \cite{kes64}, P.Walker \cite{wal77} and J.Locker \cite{lo99}.
For instance, P.Walker considered a second-order bvp with ef 
$\sin{\varrho_k}x$, such that cv
$\varrho_k:=\sqrt{\lambda_k}$ are divided into two sequences:
\[
     2{\pi}k,\ k=1,\ldots;\;\; 
    2{\pi}k + o(1),\ k=0,\ldots;\ o(1)\longrightarrow 0, \; k\longrightarrow\infty{.} 
\]
The eigenfunctions, corresponding to two \textit{close} $\varrho_k$, 
have an angle tending to zero.
Summarizing we have 
\begin{oldtheorem}[G.M.Kesel\'{}man, V.P.Mihailov, A.A.Shkalikov]
\label{thm-main-old}
\[
	L\in (SR) \Rightarrow L\in (UB), \qquad 
	L\in (WR) \Rightarrow L\in (UBP). 
\]
\end{oldtheorem}

However, further investigations failed to find even a \textit{single} bvp
with the same list of properties, if
Birkhoff-regularity is violated. Moreover, off this class the resolvent admits a polynomial and even exponential growth, see  
\cite{ef78,NaiLdoEngI} and \cite[chap.I, sec.2]{min99a}.
Of course, there is a natural candidate for \textit{good} boundary conditions: the self-adjoint ones, but they \textit{are} Birkhoff-regular \cite[$n$ even]{sal68}, 
\cite[$n$ odd]{min77a}. 

Let us recall here that \emph{essential non-selfadjointness} of bvps stems exactly
from boundary conditions and the lower order terms in $l(y)$ play a role of small perturbation.
However most of investigations of \textit{good} (in some sense) bvps deal with perturbations of the differential expression
by some subordinated functional-differential operator, see for example \cite{gr},
but not for new classes of boundary conditions,
except maybe papers \cite{yak84,yakmam,mishu97}.
\section{Green's function.}
\label{sec:green}
\subsection{Birkhoff's solutions.}
\label{sec:partsol}
Set $\varepsilon _{\jmath }=\exp (2\pi ij/n)$,
$
\varrho =\lambda ^{1/n},\; |\varrho |=|\lambda |^{1/n},
$
\begin{equation}
\arg \varrho =\arg \lambda /n,\qquad 0\leq \arg \lambda < 2\pi .
\label{eq1.6}
\end{equation}
and define sectors
\[
S_\nu=\{\varrho \,\,\,|\,\,\,\pi \nu/n\le \arg \varrho < \pi (\nu+1)/n\}.
\]
For $\varrho=\lambda^{1/n}$ we have that $\varrho \in S_0\cup S_1$. 

Let $R_0$ be a fixed positive number such that in every sector $S_\nu$ there exists a fss
 $\{y_{\jmath}(x,\varrho )\}_{j=0}^{n-1}$ of  (\ref{eq1.1a})
with an exponential asymptotics:
\begin{equation}
D^ky_{\jmath }(x,\varrho )=(\varrho\varepsilon _{\jmath })^k\cdot \exp
(i\varrho \varepsilon _{\jmath }x)[1],\ \ j,k=0,\ldots ,n-1,
 \; \abs{\varrho}\ge R_0.
  \label{eq1.7}
\end{equation}
\subsection{Canonical fss.} 
\label{sec-canon-fss}
Note that for a given sector $S_\nu$ there exists a number $p$, such that
solutions $y_{\jmath }(x,\varrho )$ decay as $j<p$ and exponentially grow
otherwise (for $x>0$), except maybe a boundary ray.
Clearly $p$ depends upon the sector's choice and values 
of $p-1$ are presented in the table \ref{table:p}.
\begin{table}
\begin{tabular}{|l|l|l|}
\hline
$\varrho \diagdown n$ & $2q$ & 2q+1 \\ \hline
$\in S_0$ & $q-1$ & $q$ \\
$\in S_1$ & $q-1$ & $q-1$ \\ \hline
\end{tabular}
\vspace{20pt}
\caption{Values of $p-1$}
\label{table:p}
\end{table}

It will be convenient to use another fss
$\{z_k\}_{k=0}^{n-1}$ 
of the equation (\ref{eq1.1a}):
\begin{equation}
   \label{eq1.10}
\boxed{
   z_k(x,\varrho ):=
      \left\{
           \begin{array}{ll}
              y_k(x,\varrho ), & \quad k=0,\ldots ,p-1, \\
              y_k(x,\varrho )/\exp (i\varrho \varepsilon _k), & \quad k=p,\ldots ,n-1.
           \end{array}
       \right. 
} 
\end{equation}
This choice is natural due to the fact
that
\begin{equation}
z_k=\mbox{O}(1)\;,k=0,\ldots ,n-1;\;\;0\le x \le 1 ,
\quad \varrho \in S_\nu.
 \label{eq1.11}
\end{equation}
\subsection{Particular solution.}
Let $W_{\jmath }(x,\varrho )$ be the algebraic complement of the element
$D^{n-1}y_{\jmath }$ in the wronskian
\[
W(x,\varrho )=\left| D^ky_{\jmath }(x,\varrho)\right| _{j,k=0}^{n-1}.
\]
Set
\(
\tilde y_{\jmath }(x,\varrho ):=W_{\jmath }/W.
\)
Calculating we find that
\begin{equation}
\tilde y_{\jmath }(x,\varrho )=\frac 1{n(\varrho \varepsilon _{\jmath
})^{n-1}}\exp (-i\varrho \varepsilon _{\jmath }x)[1].  \label{eq1.8}
\end{equation}
Introducing the kernel
\[
g_0(x,\xi ,\varrho )=i\cdot \left\{
\begin{array}{rl}
\sum\limits_{k=0}^{p-1}y_k(x,\varrho) \tilde
{y_k}(\xi,\varrho ),\ \ x>\xi &  \\
-\sum\limits_{k=p}^{n-1}y_k(x,\varrho) \tilde
{y_k}(\xi,\varrho ),\ \ x<\xi &
\end{array}
\right.
\]
we get a particular solution $g_0(f)$ of the equation $l(y)=\lambda y+f$,
\begin{equation}
g_0(f):=\int\limits_0^1g_0(x,\xi ,\varrho )f(\xi ))d\xi .  \label{eq1.9}
\end{equation}
\subsection{Formula for the Green's function.}
When applying boundary conditions \eqref{eq1.2} to a function of $x$ and $\varrho$, e.g. $z(x,\varrho)$,
it is convenient to rewrite them in vector form:
\[
   V(z) := \left(
                     \varrho^{-j} U_j(z)
             \right)_{j=0}^{n-1}.   
\]
Then define
\[
\begin{array}{lll}
g(x,\xi ,\varrho )&=& g_0(x,\xi ,\varrho )\cdot (n\varrho ^{n-1})/i, \\
H(\xi ,\varrho )   &=& V_{x}(g(x,\xi ,\varrho)).
\end{array}
\]
Here the subscript \emph{x} means that the vector boundary form $V$ acts
on the kernel $g(x,\xi ,\varrho )$ over the argument $x$.

Observe that success of Birkhoff-regularity leans heavily upon explicit formula for the 
\Grf 
\begin{align}
 \lefteqn{G(x,\xi ,\varrho )} 
 \phantom{G(x,\xi ,\varrho) } &= \frac{(-1)^n\Delta (x,\xi ,\varrho )}{n\varrho
^{n-1}\Delta(\varrho )}, 
\label{eq_green} \\
 \lefteqn{\Delta (\varrho )}  
 \phantom{G(x,\xi ,\varrho) } &= \det\gD(\varrho) =
        \det 
 \bl
        V(z_0)\ldots V(z_{n-1}) 
 \br, 
 \label{eq_chardet} \\
\Delta (x,\xi ,\varrho ) &= i\cdot \left|
\begin{array}{cc}
z^T & g(x,\xi ,\varrho ) \\
\gD(\varrho ) & H(\xi ,\varrho )
\end{array}
\right| , \label{eq1.14}
\end{align}
where $z^T$ stands for the row 
\[
   \left(z_0(x,\varrho ),\ldots ,z_{n-1}(x,\varrho )\right).
\]
$\Delta(\varrho)$ is referred to as the \emph{characteristic determinant}. 
Its estimate from below constitutes the main ingredient of the resolvent method. 
\subsection{New regularity determinants. }
\label{std_reg}
\begin{definition}
\label{def-new-reg}
Fix some $\varepsilon\in(0,\pi/{2n})$.
Let $S_\nu(\varepsilon)$ be the sector
\begin{equation}
	\label{eq_eps_sec}
	S_\nu(\varepsilon) = 
	\left\{
	       \abs{
	       \arg\varrho - \frac{(\nu+1/2)\pi}{n} 
	       }
	          \le \varepsilon 
	\right\}.
\end{equation}
Define the regularity determinants, corresponding to the sectors
$S_\nu$, via the formula
\begin{equation}
\Theta\left( S_\nu\right) := 
\lim_{\varrho \longrightarrow \infty  }\Delta(\varrho), \quad
\varrho \in S_\nu(\varepsilon).
 \label{eq_regdet}
\end{equation}
\end{definition}
Let $q=\entier{(n/2)}$. Then for $0\le k \le n-1$ set
\begin{equation}
b^i=(b_{\jmath }^i)_{j=0}^{n-1}
,\;
B_k^i=\left( b_{\jmath }^i\cdot \varepsilon _k^j\right)
_{j=0}^{n-1},\;i=0,1.  \label{eq1.3}
\end{equation}
\medskip
It is easy to calculate the limit in \eqref{eq_regdet} for the determinant
and its matrix
\begin{align}
\Theta(S_\nu) &=\Theta_p(b^0,b^1), 
    \label{eq1.4a}
    \\
 \gT(S_\nu) &=\gT_p(b^0,b^1):=
         \bl{B}_k^0,k=0,\ldots ,p-1|B_k^1,k=p,\ldots ,n-1\br.
         \label{eq1.4b}
\end{align}
The vertical line $|$ separates columns with superscripts $0$ and $1$. 
Recall that $p=p(\nu)$, $\nu=0,1$, see table \ref{table:p}.
From (\ref{eq1.4a}) it is clear that definition \ref{def-new-reg}
is equivalent to the standard one
\cite[p.361]{sal68},
but seems to be more natural and useful for generalizations.
\begin{definition}
  \label{b-reg}
We shall call boundary
conditions (\ref{eq1.2}) and the corresponding operator $L$ Birkhoff-regular
and write $L\in (R)$, if
\begin{equation}
\Theta(S_0) \ne 0,\; \Theta(S_1) \ne 0.
\label{eq1.5}
\end{equation}
\end{definition}
\begin{definition}
\label{s-reg}
Birkhoff-regular bvp is \textbf{strongly regular}, 
$L \in (SR)$, 
if
either $n$ is odd or if it is even, $n=2q$, and the second order polynomial
$F(s)$ has two simple roots, where
$F(s)=\det\F(s)$,
\[
  \F(s):=\bl B_0^0+s\cdot B_0^1,\,B_k^0,k=1,\ldots ,q-1\,|\,
  s\cdot B_q^0+B_q^1,\,
B_k^1,k=q+1,\ldots ,n-1\br
\]
Otherwise we shall call the bvp \textbf{weakly regular}
and write
$L\in (WR)$, i.e. for classes of bvps we define $(WR):=(R)\setminus(SR)$.
\end{definition}
\noindent
\subsection{Modified characteristic matrix.}
\label{subsec-mcm}
\subsubsection{Preliminaries.}
Set
\begin{align}
   \label{eq1_u_t}
\boxed{
   u_t=
   \left\{
          \begin{array}{lcll}           
          \tilde{y_t}(\xi ,\varrho ) \cdot
           n(\varrho \varepsilon _t)^{n-1} \cdot
           e^{i\varrho \varepsilon _t}
         & =  
         e^{i\varrho \varepsilon _t(1-\xi )}\cdot [1], \; & t<p\\
         \tilde {y_t}(\xi ,\varrho )\cdot n(\varrho \varepsilon _t)^{n-1} 
         & = 
         e^{i\varrho \varepsilon _t(-\xi )}\phantom{1}\cdot [1], \; & t\ge p\\
         \end{array}
   \right.  
}
\end{align}
The following formula stems immediately from definitions of $z_k$ and $u_t$:
\begin{equation}
   \label{eq-g-repr}
    g(x,\xi,\varrho) =
    \begin{cases}
       +\sum_{k<p} \varepsilon_k z_k(x,\varrho)
       u_k(\xi,\varrho)e^{-i\varrho\varepsilon_k},    
       & x > \xi, \\
       -\sum_{k\ge p} \varepsilon_k z_k(x,\varrho)
       u_k(\xi,\varrho)e^{+i\varrho\varepsilon_k},
       & x < \xi.  
    \end{cases}        
\end{equation}
Introduce notation
\[ 
   \#=\#(t)=\left\{
   \begin{tabular}{ll}
        1, & $t<p$ \\
        0, & $t\geq p$.
   \end{tabular}
\right. 
\] 
We shall omit the index if it is clear from context.
\begin{lemma}
The following representation holds true
\begin{align}
   \label{eq-H-repr}
                  H(\xi ,\varrho )  
                  &=
                  \sum_{t=0}^{n-1}(-1)^{(1-\#)} 
                  \biggl(
                         V_\#(z_t) e^{(-1)^\# i\varrho\varepsilon_t}
                  \biggr)       
                  \cdot \varepsilon_t u_t(\xi,\varrho)    \notag\\         
                  &=
                  \sum_{t=0}^{n-1}(-1)^{(1-\#)}
                  [{B}_t^\#]
                  \cdot \varepsilon_t u_t(\xi,\varrho).
\end{align}
\end{lemma}
\begin{proof}
Let $V_0$ and $V_1$ be summands of the vector functional $V$, corresponding to 
derivatives at $0$ and at $1$.
Applying $V_0$ and $V_1$ to the rhs of \eqref{eq-g-repr} over variable $x$, we get
\begin{align*}
   V_{0,x}(g) &=  - \sum_{t\ge p} \varepsilon_t u_t(\xi,\varrho)
                                V_0(z_t)e^{+i\varrho\varepsilon_t}       \\
              &=  - \sum_{t\ge p} \varepsilon_t u_t(\xi,\varrho)[B_t^0], \\           
   V_{1,x}(g) &= + \sum_{t < p}  \varepsilon_t u_t(\xi,\varrho)    
                                V_1(z_t)e^{-i\varrho\varepsilon_t}       \\
              &= + \sum_{t < p}  \varepsilon_t u_t(\xi,\varrho) [B_t^1].               
\end{align*}
Here exponentials cancel out after substituting \eqref{eq1.7}-\eqref{eq1.10}.
\end{proof}
\subsubsection{Main formula.}
\begin{lemma}
   \label{lem-green-repr}
   \begin{equation}
G(x,\xi ,\varrho )=g_0(x,\xi ,\varrho )-\frac{2\pi i}{n\varrho ^{n-1}}%
\sum_{t,k=0}^{n-1}a_{tk}(\varrho ) z_k(x,\varrho ) u_t(\xi
,\varrho ),  \label{eq1.20}
\end{equation}
where
\begin{equation}
    \label{eq1.21}
a_{tk}:=\left\{
\begin{array}{ll}
+\frac{\varepsilon _t}{2\pi }\cdot 
            \left| 
               \Delta \mathop{\longleftarrow}%
               \limits_{k} [{B}_t^1] 
            \right| 
            /\Delta , & t<p, \\
-\frac{\varepsilon _t}{2\pi }\cdot 
           \left| 
           \Delta \mathop{\longleftarrow}%
           \limits_{k} [{B}_t^0]
           \right| 
           /\Delta , & t\ge p.
\end{array}
\right.  
\end{equation}
\end{lemma}
\begin{proof}
Denote $\Delta_k$ the $k$-th column of the matrix $\gD$.
Below we use the standard agreement that $\widehat{d}$ means absence of the 
corresponding column in determinant.
Expanding $\Delta (x,\xi ,\varrho )$ along the topmost row
and replacing $H$ by the sum in \eqref{eq-H-repr},
we obtain
\begin{align*}
     & G(x,\xi,\varrho) \\
     & = 
     g_0(x,\xi,\varrho) + \frac{(-1)^n i}{n\varrho^{n-1}\Delta}                         \sum_{k=0}^{n-1}(-1)^k z_k(x,\varrho)
     \left|   
             \Delta_0\ldots\widehat{\Delta_{k}}\ldots\Delta_{n-1}H
     \right| \\
     & =
     g_0(x,\xi,\varrho) + \frac{(-1)^n i}{n\varrho^{n-1}\Delta}
     \sum_{k=0}^{n-1}(-1)^k (-1)^{n-1-k} z_k(x,\varrho)
     \left|
             \Delta\mathop{\longleftarrow}\limits_{k} H
     \right| \\
     &=
     g_0(x,\xi,\varrho) + \frac{-i}{n\varrho^{n-1}\Delta}
     \sum_{k=0}^{n-1}     
     z_k(x,\varrho)\sum_{t=0}^{n-1}
     (-1)^{(1-\#)}
     \left|
             \Delta\mathop{\longleftarrow}\limits_{k} [B_t^\#]
     \right|   \\ 
     &\phantom{g_0(x,\xi,\varrho) + \frac{-i}{n\varrho^{n-1}\Delta}}  
     \cdot\varepsilon_tu_t(\xi,\varrho),
\end{align*}
whence \eqref{eq1.20} follows.
\end{proof}
Earlier the coefficients (\ref{eq1.21})
were introduced in \cite{min99a,min95a}.
It would be quite natural to call the matrix 
\begin{equation}
\label{eq-char-matrix}
	\A=\A(\varrho)=\bl a_{tk}\br _{t,k=0}^{n-1}
\end{equation}
a \textit{\mcm} (mcm) of the
bvp (\ref{eq1.1a}),(\ref{eq1.2}) because it differs from
the analogous object in \cite[p.135]{NaiLdoEngII}
by another choice of the fss. Namely, in Naimark's
book the latter is taken analytic in $\lambda $.
\section{Further development.}
\label{sec:stone}
\subsection{Stone-regularity.}
In case of smooth coefficients of $l(y)$ the determinant $\Delta(\varrho)$ admits further terms of asymptotic expansion. Assuming that necessary amount of them is not vanishing,
investigators, A.P.Hromov, H.Benzinger, W.Eberhard, G.Freiling,  came to the notion of Stone-regularity, see for instance \cite[chap.I,sec.2]{min99a} and references therein.
This approach yields spectrum asymptotics, completeness and \textit{upper} polynomial estimate of the \Grf which is \textit{worse} than that of the Birkhoff's case. However  sharpness of this upper estimate remains open.
\subsection{Expansions of smooth functions.}
Stone-regularity continues to attract a lot of attention, see the recent books of J.Locker \cite{lo99} and R.Mennicken, M.M\"oller \cite{MM}. 
Moreover, A.A.Shkalikov and C.Tretter succeeded to establish
unconditional convergence 
for Stone-regular problems in classes of sufficiently smooth functions \cite{shk86,ShkTr}. Roughly speaking, the functions' smoothness should be enough to suppress the possible growth of the resolvent.
\subsection{Classification for $D^2$.}
For simplest two-point bvps when $l(y)\equiv{D}^2$ P.Lang 
and J.Locker \cite[p.554]{LL1} carried out 
\emph{a complete classification} of their spectral properties!
It is based on the Pl\"ucker coordinates $p_{ij}, i<j,\; i,j=0,\ldots,3$
of the matrix of coefficients in the boundary conditions (\ref{eq1.1b}).
Remind that for $n=2$ it is a $2\times 4$ matrix and 
$p_{ij}$ stands for its $2\times 2$ minor with columns $i$, $j$.
These coordinates constitute a full set of the bvp invariants and are independent up to
a well-known quadratic relation:
\[
	p_{01}p_{23}+p_{02}p_{13}+p_{03}p_{12}=0.
\]
Further in  the book \cite{lo99}
J.Locker performed a thorough investigation
of two-point bvps for $l(y)\equiv D^ny$. 
He classified degeneracy of the polynomial coefficients by the leading exponentials 
in the characteristic determinant
and obtained the same results as for classical Stone-regular bvps, also without any claim of non-spectrality.

Thus for higher order Stone-regular problems the sharpness of the resolvent's estimate in irregular cases
has not been established. Therefore it was not proved either that they
are non-spectral, not speaking of 
two-point bvps \textit{without} Stone-regularity assumptions. Such operators seem merely unattainable.

Summarizing we conclude that no \UB classification has been obtained for higher order differential operators,
cf. \cite[remark on the p.98]{lo99}.
\section{Main results.}
\label{sec:main}
\subsection{Main theorem.}
The theorem below was a widely held tacit conjecture though never formulated explicitly.
\begin{theorem}
   \label{thm-1}
   \framebox{$L \in (UB) \Rightarrow L \in (R)$}.
\end{theorem}
Theorem \ref{thm-1} together with theorem \ref{thm-main-old}
solves the \UB problem except for 
weakly regular $L$. In the latter case we have 
\UBP but need to describe the subset of \UB bvps.
Of course there are two trivial examples: 
\begin{itemize}
	\item $L\in(R)$ such that all but a finite number of ev are multiple. 
	      Then it is enough to perform Gram-Schmidt orthogonalization process
	       and get an unconditional basis. 
	\item $L$ is self-adjoint operator and $L\in (WR)$.         
\end{itemize}
\begin{problem}
Give necessary and sufficient conditions for $L \in (WR)\cap(UB)$.
\end{problem}
\subsection{Minimal resolvent's growth.}
Recall the following estimate for the resolvent
\begin{equation}
  \label{eq-lrg}
	\norm{R_\lambda} \le \frac{C}{dist(\lambda,\Lambda)} \tag{LRG}   
\end{equation}
which is called the \textsl{Linear Resolvent Growth} condition \cite{kutre01}.
It holds if $L\in (UB)$. Then inequality
\begin{equation}
\norm{R_\lambda}
\leq C\left| \varrho \right| ^{-n}
\label{eq2.1}
\end{equation}
stems from (\ref{eq-lrg}) for $\lambda$ such that 
\begin{equation}
	\label{eq_rho_dist}
	dist(\lambda,\Lambda) \ge C|\lambda|.
\end{equation}
Set 
\[
\D:=diag\left( \varepsilon _0,\ldots ,\varepsilon _{p-1},-\varepsilon
_p,\ldots ,-\varepsilon _{n-1}\right) /\left( 2\pi \right).
\]
\begin{theorem}
\label{thm-2} 
Let $L$ be a differential operator, defined by bvp (\ref{eq1.1a}),(\ref{eq1.2}). Fix $\nu \in \{ 0,1 \}$.
Given a sequence 
$
\left\{ \tau_m\right\} _1^\infty \subset
S_\nu(\varepsilon), 
$ 
such that \eqref{eq2.1} fulfills for $\varrho = \tau_m$,
we have that
\begin{align*}
	\exists \lim_{m\rightarrow \infty }\gD(\tau _m)
& = \gT_{p}\left( b^0,b^1\right), 
\\
\exists \lim_{m\rightarrow \infty }\A(\tau _m)
&=:\A_{\infty}(L)
=\gT_{p}\left( b^0,b^1\right) ^{-1}
\cdot
\gT_{p}\left(
b^1,b^0\right) \cdot \mathbf{D}
\end{align*}
and all the matrices are invertible.
\end{theorem}
Therefore in order to prove theorem \ref{thm-1}
it is enough to establish (\ref{eq2.1}) for one sequence in $S_0(\varepsilon)$,
another in $S_1(\varepsilon)$ and apply theorem \ref{thm-2}.

Converse to theorem \ref{thm-2} is also true, namely, nonvanishing of one regularity determinant yields minimal resolvent's growth in the corresponding sector 
of the $\lambda$-plane \cite{yakmam}.
\subsection{Dissipative case.}
\label{subsec:dissip}
\begin{conjecture}[S.G.Krein]
\footnote{V.A.Il'in (personal communication) has kindly informed us that it was conjectured by
S.G.Krein during one of Voronez mathematical schools in
seventies-eighties.}
Two-point dissipative bvps are Birkhoff-regular.
\end{conjecture}
\begin{theorem}
\label{thm-3}Even-order dissipative differential operators are
Birkhoff-regular.
\end{theorem}
This theorem is an immediate corollary of the theorem \ref{thm-2}.
Indeed, if $L$ is dissipative, then (\ref{eq-lrg}) is valid in the whole half-plane 
$\Cb_{-}$.
Therefore the corresponding regularity determinant 
$\varTheta(S_1)$
is nonzero. It is enough in the even-order case, since then the second determinant is the same (in our definition, see \eqref{eq1.4a} and the value of $p-1$ in the table \ref{table:p}).
\begin{remark}
\label{rem-odd-dissip}
Recently E.A.Shiryaev found another proof of theorem \ref{thm-3}
\cite{shi04}. He established self-adjointness of senior terms
of even-order dissipative boundary conditions. 
In the odd-order case he also presented an example where regularity is violated \cite[section 4]{shi04}. 
This raises a difficult question which odd order dissipative \bvps
are Birkhoff-regular (if any).

However, one nonzero determinant, namely $\varTheta(S_1)$, is enough to assert that odd-order dissipative differential operator $L_{odd}$ is \emph{half-regular} in the sense of \cite{mishu97}. 
Note that the conjugate $L_{odd}^*$ is also \emph{half-regular} by the same reasoning because its resolvent obeys 
(\ref{eq-lrg}) in $\Cb_{+}$.
\end{remark}

The following theorem sheds light on behaviour of the spectral projectors of $L_{odd}$.
\begin{theorem}
\label{thm-4}
Let $L_{odd}$ be an odd order dissipative differential operator, generated by bvp
(\ref{eq1.1a}),(\ref{eq1.2}). 
Assume that $\lambda_m$ is a simple ev (of multiplicity $1$).
For ef $u_m$ of $L_{odd}$ and the biorthogonal ef $v_m$ of $L_{odd}^*$ we have that
$(u_m,v_m)=1$.
Let $\varrho_m:=\lambda_m^{1/n}$ be the corresponding cv, $\Im\varrho_m\ge 0$.
Then the spectral projector
\[
   P_m(f) = (f,v_m)u_m,\quad f\in L^2(0,1)
\]
admits a sharp norm estimate as an operator in $L^2(0,1)$
\begin{equation}
   \label{eq-Pm-norm}
    \norm{P_m}\asymp \exp(\Im\varrho_m)/(1+\Im\varrho_m).
\end{equation}
\end{theorem}
\begin{proof}
Let $\{ \tilde{z}_k(x,\bar\varrho) \}_{k=0}^{n-1}$ be the canonical fss of the 
conjugate equation
\[
l^*(z)=\bar{\lambda}z,\; \lambda\in\Cb_{-}.
\]
Expanding $v_m$ along this fss, $u_m$ along the fss (\ref{eq1.10}) 
and invoking equation (3.3) from \cite{mishu97}, we find an asymptotic representation
\begin{align}
    u_m(x)& = c_0\cdot z_0(x,\varrho_m)[1],\\ 
    v_m(x)& = d_0\cdot \tilde{z}_0(x,\bar{\varrho_m})[1].
\end{align} 
Recall that
\begin{align}
   z_0(x,\varrho)                &= \exp(i\varrho{x})[1],\notag\\
   \tilde{z}_0(x,\bar{\varrho})  &= \exp\left(i\bar\varrho(x-1)\right)[1], \notag \\
   \norm{z_0}                    &\asymp\norm{\tilde{z_0}} 
                                    \asymp (1+\abs{\Im\varrho})^{-1/2} \notag \\
   \norm{P_m}                    &=
   \abs{c_0}\cdot\abs{d_0}\cdot\norm{z_0}\cdot\norm{\tilde{z}_0}. \notag
\end{align}
Remove for simplicity brackets since they don't affect considerations.
Then 
\[
    1=(u_m,v_m) = c_0\cdot\bar{d_0}\cdot\exp(i\varrho_m)
\]
whence estimate \eqref{eq-Pm-norm} readily follows.
\end{proof}
\section{Abstract approach.}
\label{sec:abstract}
\subsection{Functional model.}
So far we examined the resolvent approach to spectrality.
Now it's a time to look at the results for abstract linear operators 
and to try to apply them to the operator $L$.
First, remind that the spectral
theory of abstract non-selfadjoint operators has been deeply investigated, 
using M.S.Lifshi\u{c}'s characteristic function.
Presently serious attempts are made to develop 
spectral theory for
operators close to unitary \cite{kap98}, equivalently close to a self-adjoint operator, particularly translating all constructions to the language of differential operators \cite{sol89}.
However it still can not provide solution of the \UB problem.

This approach relies on the functional model theory, which is most deeply explored 
for dissipative operators \cite{SNFharm,NikShift}. 
Fundamental contribution has been done by 
A.S.Marcus, V.E.Katsnelson,   N.K.Nikolskii, B.S.Pavlov, V.I.Vasyunin and S.R.Treil, resulting finally
in a strong criterion of unconditional basicity of ef
\cite[Lect. VI, IX]{NikShift} and even of a family of invariant subspaces
\cite{tre96}. 
Let us state this remarkable result in a simple form, suitable for differential
operators.
\begin{oldtheorem}
\label{thm-dissip-um}
Assume that the differential expression (\ref{eq1.1a}) is formally self-adjoint
and $L$ is a dissipative operator. Then 
\[
L \in (UB) \ \text{in the span of eaf } \Leftrightarrow \text{uniform minimality of eaf}.
\]
\end{oldtheorem}
In abstract situation \emph{uniform minimality} is much weaker than (UB).
Nevertheless, the former seems to be unverifiable for differential operators. 
The unique class of two-point bvps
where it is known to be valid is (SR). 
However, in this case \UB is \emph{already established}!
Hence, this approach turns out to be ineffective for bvps.
\begin{remark}
From theorem \ref{thm-4} stems
\emph{uniform minimality} of ef of odd order differential operator, 
provided that cv are simple and lie in a strip:
\begin{equation}
	\label{eq_strip}
	\abs{\Im\rho_k} \le C.
\end{equation}
If we would be able to verify these assumptions, then we were able to apply theorem \ref{thm-dissip-um}. But presently neither of these conditions is possible to check.
\end{remark}
\subsection{Gubreev's development of projection method.}
Presently probably the most promising \textit{abstract} approach is worked out by G.M.Gubreev.
He succeeded to develop further Pavlov's \emph{projection method}
for finite-dimensional perturbations of Volterra dissipative operators \cite{gub00,gub03}.
To give a taste of this approach, we'll state one of his results in the simplest situation, omitting minor details for brevity.
\subsubsection{One-dimensional perturbation of $J_a$.}
Let $J_a$ be an integration operator in $L^2(0,a)$,
\[
J_{a}f(x)=\int_0^x f(t) dt, \quad f \in L^2(0,a,),
\] 
$A$ be an unbounded operator, such that
$A^{-1}$ is a one-dimensional perturbation
\[
A^{-1}h = J_ah + (h,f)g.
\]
Equivalently,
$
\mathcal{D}_A = \ker{\hat{\varphi}}
$
for some unbounded functional $\hat{\varphi}$.
Then 
\[
\varphi(z) := \hat{\varphi}\left( \exp(iaz) \right)
\]
is the generating function of the
spectrum $\Lambda$ of $A$.

Assume that $\Lambda$ lies strictly in $\Cb_{+}$, i.e obeys 
(\ref{eq-semi-bound}), and that ev are simple.
Eigenfunctions $g(\lambda_k)$ are values of 
the vector-valued entire function 
\[
g(z):=(I-zJ_a)^{-1}g.
\]
Set $w^2(x):=\| g(x) \|_{L^2(0,a)}^2$. It is a weight on $\Rb$.
Let $w_{-}(z)$ be an outer function in $\Cb_{-}$ such that
$|w_{-}(x)| = w(x),\; x\in\Rb$ \cite{garbaf}.
\begin{oldtheorem}
$A^{-1} \in (UB)$ iff
\begin{enumerate}
	\item $\varphi$ is efet with indicator diagram $[0,a]$;
	\label{it_1}
	\item $g$ coincides with restriction on $(0,a)$ of some function $g_w(x), x\in         \Rb_{+}$, see \cite{gub00}, generated by Muckenhoupt weight 
	$w^2|_{\Rb} \in (A_2)$;
	\label{it_2}
	\item $\abs{\varphi(x-i\eta)\cdot w_{-}(x-i\eta)}^2 \in (A_2),\quad$ 
	     where
	     $\eta>0$ and is fixed;
	\label{it_3}
	\item $\Lambda \in (C)$.
	\label{it_4}
\end{enumerate}
\end{oldtheorem}
Actually condition \ref{it_1} is completeness of $A$,
\ref{it_2} and \ref{it_3} are its similarity to the model dissipative operator $D_M := D$ in $L^2(\Rb)$ with domain 
\[
E=span\bigl(\{e^{i\lambda_k x}\}_{k=1}^{\infty}\bigr).
\]
Fourier transform maps $E$ onto the coinvariant subspace $H(B)\subset H_+^2$, generated by the Blaschke product 
$B(z):=\prod_{\lambda\in\Lambda}b_\lambda(z)$,
while operator $D_M$ transforms to
\[
    A_{mod} f(z) = zf(z) - \lim\limits_{z\rightarrow\infty} zf(z).
\]
This can be easily verified, applying $D_M$ to the basis elements 
$e^{i\lambda_kx}$ and performing Fourier transform.
In other words conditions 2,3 assert that $\hat\varphi$ is a right Delsarte functional
on the Sobolev space $W_2^1(0,a)$ \cite[Theorem 3.3]{gub00}.
At last, condition \ref{it_4} states that ef of $D_M$ form an unconditional basis in their span $E$.

Note that for $w(x)\equiv 1$ we have also $g\equiv 1$. 
Then $A$ can be written as
\[
    A = D, \quad \mathcal{D}_A = \ker{\hat{\varphi}} 
\]
and we return back to bases from exponentials.
\subsubsection{One-dimensional perturbation of a dissipative operator.}
Similar results are valid for one-dimensional perturbations $A^{-1}f=Bf + (f,h)g$
of abstract \underline{dissipative} Volterra operator $B$
with $(I-\lambda B)^{-1}$ being an efet. 
Define a vector-valued function 
\[
g(z):=(I-zB)^{-1}g
\]
which is called a quasi-exponential.
Then \UB is established for the family 
\[
\{g(\lambda_k)\}_{k=1}^{\infty}
\]
under assumption (\ref{eq-semi-bounda}), 
provided that weight $w^2 \in (A_2)$, where 
$w^2(s):=\| g(s)\|^2$ on $\Rb$.
The necessity of the latter condition was announced in a difficult to attain article \cite{gub94}.
Being unaware of this fact we reproved it when imaginary part of $B$ is finite-dimensional \cite{min98}. Soon after G.M.Gubreev removed this restriction \cite{gub99a}.
Moreover, condition (\ref{eq-semi-bounda}) may be weakened to the requirement that there is a horizontal strip free of spectrum \cite{gub03}
\begin{equation}
	\label{eq-free-strip}
	\abs{\Im\lambda_k}\ge h > 0.
\end{equation}
Note that for exponentials such case was treated in \cite{min92}
and served as a guideline for the proof of \UB without spectrum restrictions \cite{min91}.
\subsubsection{Finite-dimensional perturbation.}
In the review \cite{gub00} one can also find theorems, concerning \UB for linear combinations of quasi-exponentials. It is easy to see that second order bvps lead to such systems. Say, for $l(y) \equiv D^2$ the ef are
\[
   c_k g(\lambda_k) + d_k g(\lambda_k) \equiv 
   c_k \exp(i\lambda_k x) + d_k \exp(i\lambda_k x), \: x \in [0,1].
\]
Both of the systems $\{g(\lambda_k)\}_1^\infty$, $\{g(-\lambda_k)\}_1^\infty$
are candidates for \UB in $L^2(0,1)$, whence we arrive at question of building an unconditional basis from two others.
However, it is not the case when $n>2$. 
Namely, assuming for simplicity $l(y)\equiv D^n$,
we expand an ef into linear combinations of $n$ exponentials
\[
\exp(iz\varepsilon_j x), j=0,\ldots, n-1
\]
at the point $z=\lambda_k$.
Here only one ($n$ odd) or two ($n$ even) of the systems
$
     \{\exp(i\lambda_k\varepsilon_jx)\}_{k=1}^\infty
$
can be candidates for (UB). So for higher order ($n>2$) bvps the ef family no longer  fits into this scheme.

Add also that spectrum properties like $\Lambda\in (C)$ or \eqref{eq-free-strip}
with ev $\lambda_k$, replaced by cv
$\varrho_k=\lambda_k^{1/n}$,
are difficult to translate to some restrictions, imposed on boundary conditions.
Nevertheless theorem \ref{thm-A} demonstrates that for the former it is possible.
\section{Limit of mcm.}
\label{sec:mcm}
\subsection{Almost orthogonality.}
   \label{sec-almost}
An \textit{almost orthogonality} property was discovered
in \cite{min80} for ordinary differential equations.
In \cite{min93b} it was transferred to quasidifferential expressions with a summable coefficient by
the $(n-1)$-st derivative.
For Birkhoff's fss of the equation (\ref{eq1.1a}) with asymptotics (\ref{eq1.7})
\textit{almost orthogonality} asserts that
\begin{equation}  
   \label{eq-alm-ort}
\left\| \sum_{k=0}^{n-1}c_k y_k(x,\varrho ) \right\| _{L^2(0,1)}^2
\asymp \sum_{k=0}^{n-1} \abs{c_k} ^2 \norm{y_k(x,\varrho)} _{L^2(0,1)}^2
\end{equation}
for any coefficients $c_k$, which may vary with $\varrho $. 
\begin{remark}
\label{rem-almost-ut}
Moreover, the system (\ref{eq1_u_t}) is also
\emph{almost orthogonal}. This is valid because the latter has also an exponential
asymptotics and this is the unique ingredient needed for this property
\cite{min93b}. 
\end{remark}
\subsection{Boundedness of mcm}
\label{sec:bound}
Below $R_0$ is the positive number from \eqref{eq1.7}.
\begin{lemma}
\label{lem2.1}
The integral operator $g_0$ in $L^2(0,1)$, see \eqref{eq1.9}, 
admits an estimate
\begin{equation}
   \label{eq2.3}
\norm{g_0} \leq C\abs{\varrho}^{-n},\quad 
  \varrho \in S_\nu(\varepsilon),
\; \abs{\varrho} \geq R_0.
\end{equation}
\end{lemma}
\begin{proof}
Removing brackets from the asymptotic expressions for
the functions $y_k(x,\varrho )$ and $\tilde{y}_k(\xi ,\varrho )$, we
obtain a kernel $G_0(x,\xi ,\varrho )$, which naturally extends to 
$\Rb$. 
So
\[
g_0(x,\xi ,\varrho )=G_0(x,\xi ,\varrho )+O\left( \frac 1{\varrho ^n}\right).
\]
Obviously, the extended kernel coincides with the \Grf of the
self-adjoint  operator $D^n$ in $L^2(\Rb)$. 
The latter obeys an analogue of \eqref{eq2.3} in $L^2(\Rb)$.
All the more an integral operator with the kernel $G_0(x,\xi ,\varrho )$ obeys the same estimate
in $L^2(0,a)$,
which completes the proof.
\end{proof}
\begin{lemma}
\label{lem2.2}
Let $P=P(\varrho )$ be a finite dimensional operator in $%
L^2(0,1)$ with the kernel
\[
P(x,\xi ,\varrho )=\sum_{t,k=0}^{n-1}a_{tk}(\varrho )\,z_k(x,\varrho
)u_t(\xi ,\varrho ).
\]
Then 
\begin{equation}
   \label{eq2.4}
\norm{P}\leq C/\abs{\varrho},\;
\abs{\varrho}\geq R_0,\;
\varrho\in S_\nu(\varepsilon).
\end{equation}
In addition a double sided estimate holds
\begin{equation}
   \label{eq2.5}
\norm{P}\asymp \sqrt{%
\sum_{t,k=0}^{n-1}\abs{a_{tk}(\varrho )}^2\cdot \frac 1{\abs{\varrho}^2}} ,
\quad
\varrho\in S_\nu(\varepsilon),\; \abs{\varrho} \geq R_0.
\end{equation}
\end{lemma}
\begin{proof}
First observe that (\ref{eq2.4}) stems readily from (\ref{eq2.1}),
(\ref{eq2.3}) and (\ref{eq1.20}). Further, let $f\in L^2(0,1).$ Then 
\[
Pf=\sum_{k=0}^{n-1}d_k\,z_k(x,\varrho ),
\]
where
\[
d_k=\sum_{t=0}^{n-1}a_{tk}(\varrho )\int_0^1f(\xi )u_t(\xi ,\varrho )\,d\xi
.
\]
Invoking \emph{almost orthogonality} property (\ref{eq-alm-ort}),
we arrive at the relation:
\[
\left\| Pf\right\|^2_{L^2(0,1)}\asymp \sum_{k=0}^{n-1}\left| d_k\right|
^2\,\left\| z_k\right\| _{L^2(0,1)}^2.
\]
Next, a direct calculation shows that
\begin{gather}
   \label{eq-norm-z0-u0}
\norm{z_k}^2\asymp \norm{u_t}^2 \asymp \frac{1}{\abs{\varrho}},\qquad 
\varrho\in S_\nu(\varepsilon), \; \abs{\varrho}\geq R_0. 
\end{gather}
Introduce a sum
\[
   T_k = \sum_{t=0}^{n-1}\left| a_{tk}(\varrho )\right| ^2.
\]
Then \eqref{eq2.5} reduces to 
\begin{equation}
 \label{eq2.6}
\sup_{\norm{f} \leq 1}\sum_{k=0}^{n-1}\abs{d_k}^2\asymp
\frac {1}{\abs{\varrho}}\sum_{t,k=0}^{n-1}\abs{a_{tk}(\varrho )} ^2
= \frac {1}{\abs{\varrho}}\sum_{k=0}^{n-1} T_k.
\end{equation}
Fix $\varrho $ and suppose that $T_k$
attains its maximum for $k=k_0(\varrho )$. 
Then it suffices to check (\ref{eq2.6}) when its rhs reduces
to one summand:
\begin{equation}
	\label{eq-new-rhs}
	\sup_{\norm{f} \leq 1}\sum_{k=0}^{n-1}\abs{d_k}^2
	\asymp
	\frac{1}{\abs{\varrho}}T_{k_0}.
\end{equation}
Obviously the lhs of \eqref{eq2.6} has trivial bounds
\begin{equation}
\label{eq-lhs-estim}
     \sup_{\norm{f} \leq 1}\abs{d_{k_0}}^2
     \le \text{lhs \eqref{eq2.6}} \le
     \sum_{k=0}^{n-1}\sup_{\norm{f} \leq 1}\abs{d_k}^2.
\end{equation}
Applying \eqref{eq-norm-z0-u0} and remark \ref{rem-almost-ut}, we get that
for fixed $k$
\[
            \sup_{\norm{f} \leq 1}\abs{d_k}^2 
            = \norm{\sum_{t=0}^{n-1}a_{tk}\,u_t(\cdot ,\varrho )}^2_{L^2(0,1)}
            \asymp
            \sum_{t=0}^{n-1}\abs{a_{tk}}^2\norm{u_t}^2
            \asymp
            T_k\frac{1}{\abs{\varrho}}.
\] 
It allows to rewrite \eqref{eq-lhs-estim} as
\[
      \text{lhs \eqref{eq2.6}} \asymp T_{k_0}\frac{1}{\abs{\varrho}}
\]
which coincides with \eqref{eq-new-rhs}. 
\end{proof}
\begin{corollary}
The mcm is bounded
\begin{equation}
   \label{eq2.7}
\sum_{t,k=0}^{n-1}\left| a_{tk}(\varrho )\right| ^2=O\left( 1\right) ,\quad
   \varrho \in S_\nu(\varepsilon),\ \ \abs{\varrho} \geq R_0.  
\end{equation}
Indeed, one should compare (\ref{eq2.4}) and (\ref{eq2.5}).
\end{corollary}
\subsection{Proof of theorem \ref{thm-2}.}
   \label{sec:end}
Let $A_t$ be the $t$-th column of the matrix \eqref{eq-char-matrix}
\[
\A=\bl A_0,\ldots A_{n-1}\br. 
\]
In virtue of \eqref{eq1.21}  $A_t$ satisfies an equation
\begin{equation}
   \label{eq3.1}
   \gD A_t=(-1)^{(1-\#)} \cdot \frac{\varepsilon _t}{2\pi }\left[ {B}%
   _t^{\#}\right] ,\quad 0\leq t\leq n-1.
\end{equation}
\begin{lemma}
   \label{lem3.1} 
For every $t\in \left\{ 0,\ldots ,n-1\right\} $ there exists
a vector $\eta _t\in \Cb^n$ such that
\begin{equation}
\gT_p\left( b^0,b^1\right) \eta _t=
 (-1)^{(1-\#)} \cdot
 \frac{\varepsilon _t%
}{2\pi } {B}_t^{\#} .  \label{eq3.2}
\end{equation}
\end{lemma}
\begin{proof}
Fix $t\in \left\{ 0,\ldots ,n-1\right\}.$ Using compactness of the
set of vectors
\[
A_t(\varrho ),\; \varrho\in S_\nu(\varepsilon),\ \ \abs{\varrho} \geq R_0
\]
we deduce existence of a limiting vector $\eta _t$:
\begin{equation}
  \label{eq-eta-exist}
\eta _t=\lim_{l\longrightarrow \infty }A_t(\varrho_{m_l})
\end{equation}
for some subsequence $\varrho_{m_l}$. 
In the meantime the formula
\begin{equation}
  \label{eq3.3}
\lim_{l\longrightarrow \infty }\gD(\varrho_{m_l})=\gT%
_p\left( b^0,b^1\right)  
\end{equation}
stems directly from \eqref{eq1.4b}.
Combine \eqref{eq3.1},\eqref{eq-eta-exist} and \eqref{eq3.3},
 and we are done.
\end{proof}
\begin{lemma}
\label{lem3.2} Denote $R(\A)$ the image of matrix $\A$. Then
\begin{equation}
R\left( \gT\left( b^0,b^1\right) \right) \supset span\left( {B}%
_0^0,\ldots ,{B}_{p-1}^0,{B}_p^1,\ldots ,{B}_{n-1}^1\right) .  \label{eq3.4}
\end{equation}
\begin{equation}
R\left( \gT\left( b^0,b^1\right) \right) \supset span\left( {B}%
_0^1,\ldots ,{B}_{p-1}^1,{B}_p^0,\ldots ,{B}_{n-1}^0\right) .  \label{eq3.5}
\end{equation}
\end{lemma}
\begin{proof}
First, apply the matrix $\gT_p\left( b^0,b^1\right) $ to
the standard basis in $\Cb^n$ and get (\ref{eq3.4}). Second, (\ref
{eq3.5}) follows from (\ref{eq3.2}) when $t$ runs over $0,\ldots ,n-1.$ 
\end{proof}
\begin{lemma}
\label{lem3.3} The matrix $\gT_p\left( b^0,b^1\right)$ is invertible.
\end{lemma}
\begin{proof}
Set
\[
   \Q = \bl \Q^0,\Q^1\br,\;
   \Q^i=\bl B_0^i \ldots B_{n-1}^i \br,\quad 
   {\gP }=\bl \varepsilon _j^k\br _{j,k=0}^{n-1}.
\]
Inclusions (\ref{eq3.4})-(\ref{eq3.5}) yield that
\[
R\left( \gT_p\left( b^0,b^1\right) \right) \supset span \Q.
\]
But a $(j,k)$th block-entry of the product $\Q^i\cdot{\gP }^{*}$
is an $r_j\times 1$ vector
\[
b_j^i\sum_{t=0}^{n-1}\varepsilon _j^t\cdot \overline{\varepsilon _k^t}%
=b_j^i\cdot n\cdot \delta _{jk},\quad i=0,1; \quad j,k=0,\ldots ,n-1,
\]
whence
\begin{equation}
\frac{1}{n}\Q\gP^{*}=\B:=\left(
\begin{tabular}{llllll}
b$_0^0$ &  &  & b$_0^1$ &  &  \\
& $\ddots $ &  &  & $\ddots $ &  \\
&  & b$_{n-1}^0$ &  &  & b$_{n-1}^1$%
\end{tabular}
\right) .  \label{eq3.7}
\end{equation}
Due to \eqref{eq-rankB} $rank\B=\sum_{j=0}^{n-1}r_j=n$. 
Therefore 
\[
R \left( \Q\gP^{*} \right) = R(\B)=\Cb^n. 
\]
Since $\gP$ is invertible, $\Q$ 
is a full range matrix, and the same is $\gT_p\left( b^0,b^1\right)$.
\end{proof}
\section{Sparseness of cv.}
\label{sec:growth}
\subsection{Estimate off cv.}
Denote $\Gamma=\{\varrho_j\}_{j=1}^\infty$, the sequence of all distinct cv not counting multiplicities.
Fix $\nu \in \{0,1\}$ and let 
\begin{equation}
   \label{eq-g-eps}
\Gamma_\varepsilon := \Gamma \cap S_\nu(\varepsilon).
\end{equation}
Draw a hyperbolic circle
\[
    K(\varrho_j,\delta) = 
    \left\{ 
           \varrho: \abs{b_{\varrho_j}(\varrho)} \le \delta 
    \right\}
\]
around every $\varrho_j \in \Gamma$, remove them from
$S_\nu(\varepsilon)$ and denote $S_\nu(\varepsilon,\delta)$
the remaining domain.
Set 
\[
D(\varrho,\delta):=\{ \abs{z-\varrho}\le\delta\abs{\Im\varrho}\}.
\]
Below we shall often use relations from \cite[Lecture XI, formulas after (9)]{NikShift} 
\begin{align}
K(\varrho,\delta) &\supset D(\varrho,\delta), \label{eq-incl1}\\
K(\varrho,\delta) &\subset D(\varrho,\delta_1), \;
\delta_1 = \frac{2\delta}{1-\delta}. \label{eq-incl2}
\end{align}
\begin{lemma}
\label{lem-rho-est1}
Let $\varrho \in S_\nu(\varepsilon,\delta)$.
Then 
$
    \abs{\varrho -\varrho_j} \ge c\abs{\varrho},\; 
    \forall \varrho_j \in {\Gamma}.
$
\begin{proof}
Let for definiteness $\nu=0$.
If $\varrho_j\in S_1$ then 
\[
    \abs{\varrho - \varrho_j} \ge dist(\varrho,\partial S_1) 
    \ge
    \abs{\varrho}\cdot
        \sin\left( 
                   \bigl(\frac{\pi}{n}-(\frac{\pi}{2n} + \varepsilon) \bigr)
           \right) 
    \ge c\abs{\varrho}.
\]
If $\varrho_j\in S_0$
then 
$\abs{b_{\varrho}(\varrho_j)} = \abs{b_{\varrho_j}(\varrho)}>\delta$
according to the choice of $\varrho$.
So $\varrho_j\notin K_\varrho(\delta)$.
From (\ref{eq-incl1}) stems
$\varrho_j \notin D(\varrho,\delta)$,
i.e. 
\[
\abs{\varrho_j - \varrho} > \delta\abs{\Im\varrho}
\ge\delta\cdot\sin\left( \frac{\pi}{2n}-\varepsilon\right)
\] 
where we used that
$
\arg\varrho\ge \left( \frac{\pi}{2n}-\varepsilon\right).
$
\end{proof}
\end{lemma}
\begin{lemma}
\label{lem-rho-est2}
Let  $\varrho\in S_0\cup S_1$,
$\varrho_j\in\Gamma$. Then
\[
   \abs{\varrho\varepsilon_m - \varrho_j} \ge c\abs{\Im\varrho},
   \quad m=1,\ldots,n-1.
\]
\end{lemma}
\begin{proof}
Let for simplicity $m=1$. Other $m$ may be considered in a similar way.
Then $\varrho\varepsilon_1\notin S_0 \cup S_1$
while $\varrho_j$ belongs to this union. 
Assume $n>2$. Then
\[
   \abs{\varrho\varepsilon_1 - \varrho_j}
   \ge 
   dist(\varrho\varepsilon_1,\partial S_1)
   \ge
   dist(\varrho\varepsilon_1, \{\arg z=\arg \varepsilon_1\} ) = dist(\varrho,\Rb_+)
   \ge \abs{\Im\varrho}.
\]
If $n=2$ then $\varrho\varepsilon_1=-\varrho$,
$S_0\cup S_1 =\Cb_+$
and $dist(-\varrho,\Cb_+)=|\Im\varrho|$.
\end{proof}
\begin{lemma} 
\label{lem-res-est}
Let $\varrho \in S_\nu(\varepsilon,\delta)$, $\lambda=\varrho^n$. Then the estimate (\ref{eq2.1}) holds true.
\end{lemma}
\begin{proof}
We use the identity
\[
\abs{\varrho^n - \varrho_j^n}
=\prod_{m=0}^{n-1}\abs{\varrho\varepsilon_m - \varrho_j},
\]
estimate its factors for $m=0$ via lemma \ref{lem-rho-est1}, others by lemma
\ref{lem-rho-est2} and thus arrive at (\ref{eq_rho_dist}),
whence (\ref{eq2.1}) follows.
\end{proof}
\subsection{ef properties. \label{sec-ef-prop}}
Choose some circle $K(\varrho,\delta)$,
fixed an integer $N$
and take any $N$
elements from $\Gamma\cap K(\varrho,\delta)$.
Enumerate these cv $\varrho_1,\ldots,\varrho_N$
and denote respective ef $u_1,\ldots,u_N$.
Set
\begin{equation}
	\label{eq_p2}	
	\omega_{lq}(x,\varrho) := \frac{1}{q!}\frac{d^q}{d\varrho^q}\omega_{l0}(x,\varrho),
	\; l=0,\ldots,n-1	
\end{equation}
where $\omega_{l0}(x,\varrho) := z_l(x,\varrho)$. 
We used functions $\omega_{lq}$ extensively in 
\cite[chap.4]{min99a}.
It is easy to verify the estimate \cite[chap.4, lemma 5.1]{min99a}
\begin{equation}
\label{eq_p3}
	\norm{ \omega_{lq}(x,\varrho) }^2 \asymp \frac{1}{\abs{\varrho}^{2q+1}}, \ \varrho\in S_\nu(\varepsilon).
\end{equation}
Next, expanding the ef $u_j(x)$ over the system (\ref{eq1.10}), we get
\begin{equation}
	\label{eq_p4}
	u_j(x) = \sum_{l=0}^{n-1} d_{jl} \omega_{l0}(x,\varrho_j),\ j=1,\ldots,N.
\end{equation}
Due to \textit{almost orthogonality} of the fss (\ref{eq1.10})
\begin{equation}
	\label{eq_p5}
	\norm{ u_j }^2 \asymp \sum_{l=0}^{n-1} \abs{d_{jl}}^2\cdot\norm{\omega_{l0}}^2 
	            \asymp \sum_{l=0}^{n-1} \abs{d_{jl}}^2\cdot \frac{1}{\abs{\varrho_j}}.
\end{equation}
\subsection{Norm of a linear combination of ef.}
Now let $u(x)$ be some linear combination of the ef $\{u_j(x)\}_1^N$
\begin{equation}
	\label{eq_p6}
	u(x) = \sum_{j=1}^N c_j\cdot u_j(x)
\end{equation}
Normalize the ef $u_j,j=1,\ldots,N$.
Since they are subsystem of an unconditional basis, 
the norm equivalence holds true
\begin{equation}
	\label{eq_p7}
	\norm{u}^2 \asymp \sum_{j=1}^N \abs{c_j}^2\cdot \norm{u_j}^2
	\asymp 	\sum_{j=1}^N \abs{c_j}^2.
\end{equation}
\subsection{Canonical ef representation.}
From the other hand, the following representation is also valid  \cite[chap.4]{min99a}
\begin{equation}
	\label{eq_p8}
	u(x) = \sum_{l=0}^{n-1}\sum_{m=0}^{\infty} a_{ml} \cdot \omega_{lm}(x,\varrho)
\end{equation}
where
\begin{equation}
	\label{eq_p9}
	a_{ml} = \sum_{j=1}^N c_j\cdot d_{jl}\cdot (\varrho_j - \varrho)^m.
\end{equation}
We merely substituted (\ref{eq_p4}) into (\ref{eq_p6}) and
expanded every summand $\omega_{l0}(x,\varrho_j)$ into Taylor series centered in $\varrho$. Moreover, 
the following estimate is valid \cite[chap.4]{min99a}
\begin{equation}
	\label{eq_p10}
	\norm{u}^2 \asymp \sum_{l=0}^{n-1}\sum_{m=0}^{N-1} \abs{a_{ml}}^2 \cdot 
	        \norm{\omega_{lm}(x,\varrho)}^2
	        \asymp \sum_{l=0}^{n-1}\sum_{m=0}^{N-1} \abs{a_{ml}}^2 \cdot      \frac{1}{\abs{\varrho}^{2m+1}}, \;\varrho\in\Gamma_\varepsilon{.}
\end{equation}
In fact 
there we assumed that $\varrho$ lies in a strip 
$\abs{\Im\varrho} \le C$
but the case 
$\varrho\in S_\nu(\varepsilon)$ may be considered along the same lines and even simpler.
Note that the row $(a_{00},\ldots,a_{0,n-1})$ coincides
with linear combination of matrix 
$d$ rows, 
\begin{equation}
\label{eq-matrix-d}
d=\left[ d_{jl}\right]_{j=1,l=0}^{\ N\;\, n-1}.
\end{equation}
\subsection{Spectrum sparseness.}
Recall a definition of a sparse sequence.
\begin{definition}
Let $P$ be a sequence of points in $\Cb_{+}$. Then
  $P \in (S)$ if for some $\delta > 0$
\[
   K(\varrho,\delta) \cap K(\mu,\delta) = \emptyset,
   \; \;
   \varrho \neq \mu, \  \varrho,\mu \in P.
\]  
Equivalently each circle $K(\varrho,\delta)$ contains $\le 1$ element from
$P$.
\end{definition}
\begin{definition}
\label{def-nS}
$P$ is an $N$-sparse sequence,
$P \in (NS)$ if for some $\delta > 0$
\[
   \#\left(P \cap K(\varrho,\delta)\right) \le N, \; \forall\varrho\in S_0\cup S_1.
\]
\end{definition}
\begin{lemma}
   \label{lem-nS}
Fix $\nu\in\{0,1\}$ and consider the set
$\Gamma_\varepsilon$, see (\ref{eq-g-eps}).
Then for sufficiently small $\delta$ in definition \ref{def-nS}
$\Gamma_\varepsilon \in (nS)$.
\end{lemma}
\begin{proof}
Choose some circle $K(\varrho,\delta)$ intersecting $\Gamma_\varepsilon$.
Assume on the contrary that 
$\#\left( K(\varrho,\delta) \cap\Gamma_\varepsilon\right)>n$.
Take any $N=n+1$ cv from this intersection
and enumerate them
$\{\varrho_j\}_{j=1}^N$. 
They are distinct because we chose $\Gamma$ not counting multiplicities.
Then the rows of the matrix (\ref{eq-matrix-d}) are linearly dependent.
Therefore for appropriate coefficients $c_j$
in (\ref{eq_p9}) we have 
\begin{equation}
\label{eq_pa}
	a_{0l}=0,\ l=0,\ldots,n-1.
\end{equation}
Observe that for any $\mu$ in the larger circle $D(\varrho,\delta_1)$ $\max$ and $\min$
of the ratio $\abs{\mu}/\abs{\varrho}$ are attained when $\abs{\mu}=\abs{\varrho}\pm\delta_1\abs{\Im\varrho}$.
Therefore $\abs{\mu}/\abs{\varrho}\in [1-\delta_1,1+\delta_1]$ whence
\begin{equation}
	\label{eq-rho_j-equiv}
	\abs{\varrho_j} \asymp \abs{\varrho}, \; j=1,\ldots,N.
\end{equation}
Normalize $c_j$
\begin{equation}
	\label{eq_p11}
	\sum_{j=1}^{N} \abs{c_j}^2 = 1.
\end{equation}
Then from (\ref{eq_p9}),(\ref{eq_p11}) and (\ref{eq-incl2})
stems an estimate
\begin{equation}
	\label{eq_p12}
	\abs{a_{ml}}^2 \le \sum_{j=1}^{N} \abs{d_{jl}}^2\cdot         \abs{\delta_1\cdot\Im\varrho}^{2m}\le
	\sum_{j=1}^{N} \abs{d_{jl}}^2\cdot\delta_1^{2m}
	\cdot\abs{\varrho}^{2m}.
\end{equation}
Substituting (\ref{eq_p12}) into (\ref{eq_p10}) and taking into account (\ref{eq_pa}),
(\ref{eq-rho_j-equiv}) and (\ref{eq_p5}),
we obtain
\begin{align}
	\label{eq_pb}
  \norm{u}^2 & \le C\cdot 
	\sum_{l=0}^{n-1} \sum_{m=1}^{N-1} \sum_{j=1}^N
	\abs{d_{jl}}^2 \cdot \delta_1^{2m} \cdot 
	\frac{\abs{\varrho}^{2m}}{\abs{\varrho}^{2m+1}} \notag \\
	& \le  
	C\cdot \delta_1^2 \cdot \sum_{l=0}^{n-1} \sum_{j=1}^N
	\abs{d_{jl}}^2 \cdot \frac{1}{\abs{\varrho}} \notag\\
	& \le
	C\cdot \delta_1^2 \cdot \sum_{j=1}^N \norm{u_j}^2
	= (n+1)C\delta_1^2.
\end{align}
The latter contradicts (\ref{eq_p7}), (\ref{eq_p11})
provided $\delta_1$ (from (\ref{eq-incl2}) ) is sufficiently small.
\end{proof}
\section{\UB conjecture.}
\label{sec:proof-thm1}
For the sake of definiteness let $\nu=0$. The case $\nu=1$ may be considered similarly.
Define 
\[
\mathcal{D}=S_0\bigl( \varepsilon \bigr) \cap \{ r\le \abs{\varrho} \le r+{\delta}r \}.
\] 
Recall that the boundary of $S_0\bigl(\varepsilon\bigr)$ are the rays $\arg\varrho=\frac{\pi}{2n}\pm\varepsilon$.
Set $\omega=(2\varepsilon)/M$. We will choose the integer $M$ later.
Dissect $S_0\bigl(\varepsilon\bigr)$ by $M$ rays 
\[
T_k=\{ \arg\varrho=k\omega + \frac{\pi}{2n}-\varepsilon\}, \; k=0,\ldots,M.
\]
They divide $\mathcal{D}$ into $M$ quadrilaterals $\Pi_k,\; k=0,\ldots,M-1$.
\subsection{Quadrilaterals.}
\begin{lemma} 
\label{lem-4.1}
$diam \Pi \le 2\delta r$.
\begin{proof}
Since $\Pi_k$ lies in an angle of the opening $\omega$
and has vertexes
\begin{align*}
V_1 &= r\cdot \exp(i\theta),
V_2 = r\cdot \exp(i(\theta+\omega)),
V_3 = (r+\delta r)\cdot \exp(i(\theta+\omega)),\\
V_4 &= (r+\delta r)\cdot \exp(i\theta), \quad
\theta = \frac{\pi}{2n}-\varepsilon+\omega{k},
\end{align*}
then its diameter coincides with $\abs{V_3-V_1}=\abs{V_4-V_2}$.
Denote $\stackrel{\frown}{V_1 V_2}$ the arc between the points $V_1, V_2$. Let module 
$\abs{\ }$ stand also for the arc's length. Then 
\[
   \abs{V_3 - V_1} \le \abs{V_3-V_2}+\abs{\stackrel{\frown}{V_1 V_2}} \le \delta{r} + \omega{r} \le 2\delta{r},
\]
if we take $M\ge \frac{2\varepsilon}{\delta}$, say
$M=\entier{(\frac{2\varepsilon}{\delta}})+1$.
\end{proof}
\end{lemma}
Let $Q$ be the center of the interval $[V_1,V_3]$.
Then the circle $D(Q,\delta_2)$
contains quadrilateral $\Pi_k$ for appropriate $\delta_2$.
Namely, it is enough if its radius $r_Q \ge 2\delta{r}$. But
\[
   r_Q = \delta_2\cdot\Im{Q} \ge \delta_2\cdot{r}\sin(\frac{\pi}{2n}-\varepsilon).
\]
So we can take $\delta_2 = 2\delta/\sin(\frac{\pi}{2n}-\varepsilon)$.
\begin{lemma}
\label{lem-4.2}
Number of cv in $\mathcal{D} < \frac{n}{\delta}$.
\begin{proof}
In angular direction $\mathcal{D}$ is covered by $M$ angles of opening $\omega$.
An intersection of $\mathcal{D}$ with each angle lies in $D(Q,\delta_2)$.
Observe that in turn this circle is contained in $K(Q,\delta_2)$, see (\ref{eq-incl1}).
Choose $\delta_2$ as is needed in lemma \ref{lem-nS}. Then, according to this lemma
$K(Q,\delta_2)$ contains no more than $n$ cv from $\Gamma_\varepsilon$. 
Hence, assuming $2\varepsilon+\delta<1$
\[
\#(\Gamma\cap\mathcal{D})
\le Mn \le (\frac{2\varepsilon}{\delta}+1)\cdot n < \frac{n}{\delta}.
\]
\end{proof}
\end{lemma}
\subsection{Areas' estimates.}
Reenumerate the cv in $\mathcal{D}$:
\[
   \varrho_1,\ldots,\varrho_N,\; N \le N_0:=\entier{(\frac{n}{\delta})}+1
\]
and set
\[
    K = \bigcup_{j=1}^{N} D(\varrho_j,\frac{\delta}{N_0}).
\]
\begin{lemma}
\label{lem-4.3}
A relative area's estimate holds
\begin{equation}
\label{eq-4.1}
    \abs{K}/\abs{\mathcal{D}} \le C\delta^2/\varepsilon{.}
\end{equation}
\begin{proof}
From one hand the area of $\mathcal{D}$ is subject to an estimate
\[
  \abs{\mathcal{D}}=
   \varepsilon \cdot
  \bigl( r^2(1+\delta)^2-r^2 \bigr) 
  > \varepsilon\cdot 2r^2\delta.
\]
From the other hand the area of $K$ admits an upper estimate:
\[
\begin{array}{lllll}
 \abs{K} \le &  \sum_{j=1}^N \abs{D(\varrho_j,\frac{\delta}{N_0})} 
  & = &  \sum_{j=1}^N \pi\cdot\bigl| \frac{\delta}{N_0} \Im\varrho_j \bigr|^2 
  \\  
  \phantom{|K|} \le & N\pi\frac{\delta^2}{N_0^2}\cdot\max\abs{\varrho_j}^2 
  & \le &
   N\pi\frac{\delta^2}{N_0^2}\cdot \bigl( r(1+\delta) \bigr)^2 
   \\   
  \phantom{|K|} \le &
   \pi\delta^2(1+\delta)^2
   r^2/N_0   
   &&
\end{array}
\]
whence
\[
   \abs{K}/\abs{\mathcal{D}} \le 
   \frac{\pi\delta^2(1+\delta)^2r^2/N_0}
   {\varepsilon\cdot2\delta r^2} 
   \le C\delta/(N_0\varepsilon) 
   \le C\delta \frac{\delta}{n\varepsilon} 
   = C\delta^2/{\varepsilon}.
\]
\end{proof}
\end{lemma}
We have also to take into account the area of the circles $K(\varrho_j,\delta)$,
intersecting $\mathcal{D}$, such that their \textit{centers} $\varrho_j$
lie outside $S_0(\varepsilon)$.
Replacing them with greater ones $D(\varrho_j,\delta_1)$, we see that their area
attains maximum if their \textit{centers} lie on the boundary rays of $S_0(\varepsilon)$. Consider for instance one of them, namely the ray
${Ray}_1 =\{ \arg\varrho=\frac{\pi}{2n}+\varepsilon\}$.
Then the corresponding intersections are inside the angle 
\[
Ang_1=\{  \frac{\pi}{2n}+\varepsilon -\nu \le \arg{z} \le \frac{\pi}{2n}+\varepsilon \}. 
\]
\begin{lemma}
\label{lem-4.4} The opening 
$\nu\le c\delta_1$.
\begin{proof}
Take $\varrho\in{Ray}_1$ and draw a circle 
$D(\varrho,\delta_1)$.
Clearly ${Ray}_1$ is its tangent at some point $A$.
So the radius 
$\overrightarrow{\varrho{A}}$ 
is perpendicular to the ${Ray}_1$
at $A$. Thus
\[
\delta_1\abs{\Im\varrho} = \abs{\overrightarrow{\varrho{A}}} = 
\abs{\varrho}\cdot\sin{\nu}
\]
whence
\[
   \sin\nu = \delta_1\frac{\abs{\Im\varrho}}{\abs{\varrho}}
           = \delta_1\sin(\frac{\pi}{2n}+\varepsilon)
\]
and we are done since 
$\sin\nu\ge \frac{2}{\pi}\nu, \ \nu \in (0,\pi/2)$.
\end{proof}
\end{lemma}
The intersection with circles $D(\varrho,\delta_1)$
near the other boundary ray of $S_0(\varepsilon)$
is estimated along the same lines and is contained
within the angle 
\[
Ang_0=\{
\frac{\pi}{2n} -\varepsilon 
\le \arg z \le 
\frac{\pi}{2n} -\varepsilon + \nu\}.
\]
\subsection{Notemptiness of \emph{good} domains.}
Set $Ang=Ang_0\cup Ang_1$.
\begin{lemma}
\label{lem-4.5} The \emph{good} domain 
$\mathcal{D}\setminus K$ is not empty.
\begin{proof}
First note that
\begin{equation}
\label{eq-4.2}
\abs{\mathcal{D}\cap Ang}/\abs{\mathcal{D}}
\le
  \frac{2\nu}{2\varepsilon}
\le \frac{c\delta_1}{\varepsilon}.
\end{equation}
Summing up the rhs of inequalities 
(\ref{eq-4.1})-(\ref{eq-4.2}) we get
\[
C\delta^2/\varepsilon+c\delta_1/\varepsilon = 
C\delta^2/\varepsilon+c\frac{2\delta}{(1-\delta)\varepsilon}.
\]
Since $\varepsilon$ is fixed, this expression can be made as small as needed
for sufficiently small $\delta$.
\end{proof}
\end{lemma}
\subsection{Completion of the proof.}
At last, setting $r_m=(1+\delta)^m,\; m=1,2,\ldots$,
we get a sequence of domains
$\mathcal{D}_m$ and thus a sequence of points
$\tau_m \in \mathcal{D}_m \subset S_\nu(\varepsilon,\delta)$, tending to infinity. 
According to lemma \ref{lem-res-est}
the resolvent obeys the estimate \eqref{eq2.1} with $\varrho=\tau_m$.
It suffices to apply theorem \ref{thm-2}
and theorem \ref{thm-1} is proved.

\textbf{Acknowledgement.}
We take an opportunity to thank
B.S.Pavlov, G.M.Gubreev, S.A.Avdonin and S.A.Ivanov for providing reprints and 
Yu.Lyubarskii, G.Freiling and M.M\"oller for clarifications.

\def\cprime{$'$}
\providecommand{\bysame}{\leavevmode\hbox to3em{\hrulefill}\thinspace}
\providecommand{\MR}{\relax\ifhmode\unskip\space\fi MR }
\providecommand{\MRhref}[2]{%
  \href{http://www.ams.org/mathscinet-getitem?mr=#1}{#2}
}
\providecommand{\href}[2]{#2}

\end{document}